June 15, 2022 update

# Bernoulli, poly-Bernoulli, and Cauchy polynomials in terms of Stirling and r-Stirling numbers

**Khristo N. Boyadzhiev**
Department of Mathematics
Ohio Northern University, Ada, OH 45810, USA
k-boyadzhiev@onu.edu

*Abstract.* We review and discuss some results on the representation of Bernoulli, poly-Bernoulli numbers, and Bernoulli and Cauchy polynomials in terms of Stirling numbers of the first or second kind, or in terms of r-Stirling numbers.

## 1. Introduction

More than forty years ago Henry Gould published a very informative paper [7] discussing the classical formula

(1) $$B_n = \sum_{k=0}^{n} \frac{1}{k+1} \sum_{j=0}^{k} \binom{k}{j} (-1)^j j^n ,$$

where $B_n$ are the Bernoulli numbers. Among other things, Gould indicated that the formula was very old and was known at least in 1893. It appears also in Hasse [11]

The representation can be written as

(2) $$B_n = \sum_{k=0}^{n} \frac{(-1)^k k!}{k+1} S(n,k)$$

where

$$S(n,k)=\frac{(-1)^k}{k!}\sum_{j=0}^{k}\binom{k}{j}(-1)^j j^n$$

are the Stirling numbers of the second kind, originating in the works of James Stirling (1693-1770) (see [2] for the history of these numbers and for a short proof of (2)). Further details and combinatorial interpretation can be found in [5, 8].

It would be nice to extend this result to the Bernoulli polynomials

$$B_n(x)=\sum_{k=0}^{n}\binom{n}{p}B_p x^{n-p} \quad .$$

Indeed, if we replace here $B_p$ by their representations (2) and change the order of summation we find

$$B_n(x)=\sum_{k=0}^{n}\frac{(-1)^k k!}{k+1}\left\{\sum_{p=0}^{n}\binom{n}{p}S(p,k)x^{n-p}\right\} \tag{3}$$

which directly extends (2), as $B_n(0)=B_n$.

Another extension of (2) was obtained recently by Guo, Mezo, and Qi [10] (see also Neto [17]). Namely, Guo et al found the representation

$$B_n(r)=\sum_{k=0}^{n}\frac{(-1)^k k!}{k+1}S_r(n+r,k+r) \tag{4}$$

for all integers $n,r\geq 0$, where $S_r(n,m)$ are the r-Stirling numbers, extending $S(n,m)$. In combinatorics, $S_r(n,m)$ is the number of ways to partition the set $\{1,2,...,n\}$ into $m$ nonempty disjoint subsets so that the numbers $1,2,...,r$ are in different subsets. Thus $S_0(n,m)=S(n,m)$, and when $r=0$, the representation (4) turns into (2).

The r-Stirling numbers were studied in detail by Broder in [4]. Equation (32) in Broder's paper [4] gives the following representation of r-Stirling numbers in terms of Stirling numbers of the second kind

$$S_r(n+r,k+r)=\sum_{p=0}^{n}\binom{n}{p}S(p,k)\,r^{n-p} \,. \tag{5}$$

Here in the right hand side we recognize the sum in the braces in (3). Thus (4) follows immediately from (3) and Broder's equation (5).

In the next section we discuss another result, possibly originating from Niels Nielsen [18], which extends directly representation (1). In the third section we extend the representation to poly-Bernoulli polynomials. In section 4 we give some brief notes on Cauchy numbers and polynomials.

It is good to mention that together with (2) there is another expression of the Bernoulli numbers in terms of $S(n,k)$

$$B_n = \frac{n}{2(2^n - 1)} \sum_{k=0}^{n-1} \frac{(-1)^k k!}{2^k} S(n, k+1) .$$

It appeared in Garabedian's paper [6] (see also [7]).

**Remark**. Let $s(n,k)$ be the Stirling numbers of the first kind [5, 8] (see also Section 4). It is known that for every $n, p$

$$\sum_{k=0}^{n} S(n,k) s(k,p) = \delta_{n,p} \quad \text{and} \quad \sum_{k=0}^{n} s(n,k) S(k,p) = \delta_{n,p} .$$

Because of this duality, if $a_n$, $b_n$ $(n = 0,1,...)$ are sequences, we have

$$b_n = \sum_{k=0}^{n} S(n,k)\, a_k \;\Leftrightarrow\; a_n = \sum_{k=0}^{n} s(n,k)\, b_k$$

(Stirling transform and its inverse). It follows from (2) that

$$\sum_{k=0}^{n} s(n,k) B_k = \frac{(-1)^n n!}{n+1} = D_n .$$

The numbers $D_n$ are known as the Daehee numbers. This remark applies also to equation (7) below and some others.

## 2. Nielsen's representation

Todorov [19] obtained the formula

$$B_n(x) = \sum_{k=0}^{n} \frac{(-1)^k \Delta^k x^n}{k+1}$$

where $\Delta f(x) = f(x+1) - (x)$ is the finite difference. It is well known that for any function $f(x)$,

$$\Delta^k f(x) = \sum_{j=0}^{k} \binom{k}{j} (-1)^{k-j} f(x+j)$$

and thus Todorov's result can be written in the form

**Proposition 1**. *For every* $n \geq 0$ *and every* $x$*, the Bernoulli polynomials have the representation*

(6) $$B_n(x) = \sum_{k=0}^{n} \frac{1}{k+1} \left\{ \sum_{j=0}^{k} \binom{k}{j} (-1)^j (x+j)^n \right\} .$$

This is a direct extension of (1), as $B_n(0) = B_n$. The representation was also obtained by Guillera and Sondow in [9] by extending the work of Hasse [11]. A recent independent proof was given in [3]. The result, however, is much older; it can be recognized in equation (18) on page 232 in Nielsen's book [18].

Proposition 1 implies the representation (3) and therefore, extends also the Guo-Mezo-Qi result. We have:

**Lemma 2**. *For every* $n \geq 0$ *and* $0 \leq k \leq n$*,*

$$\sum_{j=0}^{k} \binom{k}{j} (-1)^j (x+j)^n = (-1)^k k! \sum_{p=0}^{n} \binom{n}{p} S(p,k) x^{n-p}$$

*and when* $x = r$ *is a non-negative integer,*

$$\sum_{j=0}^{k} \binom{k}{j} (-1)^j (r+j)^n = (-1)^k k! S_r(n+r, k+r) .$$

The proof is trivial. By expanding the binomial and then changing the order of summation we find

$$\sum_{j=0}^{k} \binom{k}{j} (-1)^j (x+j)^n = \sum_{j=0}^{k} \binom{k}{j} (-1)^j \sum_{p=0}^{n} \binom{n}{p} x^{n-p} j^p$$

$$= \sum_{p=0}^{n} \binom{n}{p} x^{n-p} \left\{ \sum_{j=0}^{k} \binom{k}{j} (-1)^j j^p \right\} = \sum_{p=0}^{n} \binom{n}{p} x^{n-p} \left\{ k!(-1)^k S(p,k) \right\} .$$

The second equation comes from Broder's formula (5).

## 3. Poly-Bernoulli numbers and polynomials

The poly-Bernoulli numbers and polynomials were introduced by Kaneko [12]. The poly-Bernoulli polynomials $\mathrm{B}_m^{(q)}(x)$, $q \ge 1$, can be defined by the generating function

$$\frac{\mathrm{Li}_q(1-e^{-t})}{1-e^{-t}}\, e^{xt} = \sum_{n=0}^{\infty} \mathrm{B}_n^{(q)}(x) \frac{t^n}{n!}$$

(see Bayad and Hamahata [1]). Here $\mathrm{Li}_q(x) = \sum_{n=1}^{\infty} \frac{x^n}{n^q}$ is the polylogarithmic function. When $q=1$, we have $\mathrm{Li}_1(1-e^{-t}) = t$ and $\mathrm{B}_n^{(1)}(x) = (-1)^n B_n(-x)$. The numbers $\mathrm{B}_n^{(q)} = \mathrm{B}_n^{(q)}(0)$ are the poly-Bernoulli numbers. Clearly $\mathrm{B}_n^{(1)} = (-1)^n \mathrm{B}_n$

Kaneko showed that the poly-Bernoulli numbers can be written in terms of Stirling numbers of the second kind,

$$\mathrm{B}_n^{(q)} = (-1)^n \sum_{k=0}^{n} \frac{(-1)^k k!}{(k+1)^q} S(n,k)\,. \tag{7}$$

This representation can be extended to poly-Bernoulli polynomials by the same method as above. Conveniently, Bayad and Hamahata [1] showed that

$$\mathrm{B}_n^{(q)}(x) = \sum_{k=0}^{n} \frac{1}{(k+1)^q} \sum_{j=0}^{k} \binom{k}{j} (-1)^j (x-j)^n$$

(for another proof see [3]). Expanding the binomial and changing the order of summation we find

$$\sum_{j=0}^{k} \binom{k}{j} (-1)^j (x-j)^n = \sum_{j=0}^{k} \binom{k}{j} (-1)^j \sum_{p=0}^{n} \binom{n}{p} (-1)^p x^{n-p} j^p$$

$$= \sum_{p=0}^{n} \binom{n}{p} (-1)^p x^{n-p} \left\{ \sum_{j=0}^{k} \binom{k}{j} (-1)^j j^p \right\} = \sum_{p=0}^{n} \binom{n}{p} (-1)^p x^{n-p} \left\{ k!(-1)^k S(p,k) \right\},$$

and this brings to the next result:

**Proposition 4.** *For any* $n \ge 0$, $q \ge 1$,

$$\mathrm{B}_n^{(q)}(x) = \sum_{k=0}^{n} \frac{k!(-1)^k}{(k+1)^q} \left\{ \sum_{p=0}^{n} \binom{n}{p} S(p,k)(-1)^p x^{n-p} \right\}.$$

This representation turns into (7) for $x = 0$ and into (3) for $q = 1$.

Substituting here $x = -r$ , $r \geq 0$ an integer, we obtain the poly-Bernoulli analog of formula (4), expressing poly-Bermoulli polynomials in terms of r-Stirling numbers.

**Corollary 5.** *For any integers* $n \geq 0$, $r \geq 0$, $q \geq 1$, *we have*

$$\mathrm{B}_n^{(q)}(-r) = (-1)^n \sum_{k=0}^{n} \frac{k!(-1)^k}{(k+1)^q} \left\{ \sum_{p=0}^{n} \binom{n}{p} S(p,k)\, r^{n-p} \right\}$$

$$= (-1)^n \sum_{k=0}^{n} \frac{k!(-1)^k}{(k+1)^q} S_r(n+r, k+r).$$

## 4. Cauchy numbers and polynomials

We mention here very briefly some similar developments with Cauchy numbers and polynomials. The Cauchy numbers of the first kind $c_n$, $n = 0,1,...$ , are defined by the exponential generating function

$$\frac{x}{\ln(1+x)} = \sum_{n=0}^{\infty} c_n \frac{x^n}{n!}$$

or directly by the integral

$$c_n = \int_0^1 x(x-1)...(x-n+1)\, dx \ .$$

Recalling that the Stirling numbers of the first kind $s(n,k)$ are defined by the generating function (see [5], [8])

$$x(x-1)...(x-n+1) = \sum_{k=0}^{n} s(n,k) x^k$$

we find immediately the representation

$$c_n = \sum_{k=0}^{n} \frac{s(n,k)}{k+1}$$

which somewhat resembles (2). Information about Cauchy numbers can be found in [5], [14], [15], and [16]. For instance, the above representation is given on p. 294 of [5]. Komatsu [14] has defined poly-Cauchy numbers of the first kind $c_n^{(q)}$, $q \geq 1$, which have the representation

$$c_n^{(q)} = \sum_{k=0}^{n} \frac{s(n,k)}{(k+1)^q}$$

similar to (7). The Cauchy polynomials $c_n(z)$ of the first kind are defined by the function

$$\frac{x}{(1+x)^z \ln(1+x)} = \sum_{n=0}^{\infty} c_n(z) \frac{x^n}{n!},$$

and recently Komatsu and Mezo [15] obtained a representation of these polynomials in terms of r-Stirling numbers of the first kind $s_r(n,k)$

$$c_n(r) = \sum_{k=0}^{n} \frac{s_r(n+r,k+r)}{k+1} \quad .$$

Here $s_r(n,k) = (-1)^{n-k} \begin{bmatrix} n \\ k \end{bmatrix}_r$, where the unsigned r-Stiring numbers of the first kind $\begin{bmatrix} n \\ k \end{bmatrix}_r$ are described in [15]. This representation resembles (4). More details and other similar results can be found in [15].

## 5. Higher order Bernoulli numbers and Stirling numbers

The higher order Bernoulli numbers $B_n^{(k)}$, $n \geq 0, k \geq 1$ (not to be confused with the poly-Bernoulli numbers) are defined by the generating function

$$\left(\frac{t}{e^t-1}\right)^k = \sum_{n=0}^{\infty} B_m^{(k)} \frac{t^n}{n!}.$$

For $k=1$ these are the usual Bernoulli numbers. The following representation was obtained by Kim et al. [13].

**Proposition 6.** *For* $m \geq 0,\ k \geq 1$ *we have*

$$B_m^{(k)} = \sum_{n=0}^{m} s(m+k,k)\binom{m+k}{k}^{-1} S(m,n)\,. \tag{8}$$

When $k=1$ this turns into the classical representation (2) as $s(m+1,1)=(-1)^m m!$.

*Proof.* Recall the exponential generating functions of the Stirling numbers of the first and second kind

$$\frac{(\ln(1+t))^k}{k!} = \sum_{n=0}^{\infty} s(n,k)\frac{t^n}{n!},$$

$$\frac{(e^t-1)^k}{k!} = \sum_{n=0}^{\infty} S(n,k)\frac{t^n}{n!}\ .$$

We have (with $n-k=m$ in the second equality)

$$\left(\frac{\ln(1+t)}{t}\right)^k = \sum_{n=k}^{\infty} s(n,k)k!\frac{t^{n-k}}{n!} = \sum_{m=0}^{\infty} s(m+k,k)\frac{k!}{(m+k)!}t^m = \sum_{m=0}^{\infty} s(m+k,k)\binom{m+k}{k}^{-1}\frac{t^m}{m!}\ .$$

The numbers

$$D_m^{(k)} = s(m+k,k)\binom{m+k}{k}^{-1}$$

are known as the higher order Daehee numbers Thus we have

$$\left(\frac{\ln(1+t)}{t}\right)^k = \sum_{n=0}^{\infty} D_n^{(k)} \frac{t^n}{n!}$$

Replacing here $t$ by $e^t-1$ we find

$$\sum_{n=0}^{\infty} D_n^{(k)} \frac{(e^t-1)^n}{n!} = \left(\frac{t}{e^t-1}\right)^k = \sum_{m=0}^{\infty} B_m^{(k)} \frac{t^m}{m!}$$

and at the same time (changing the order of summation in the second term)

$$\sum_{n=0}^{\infty} D_n^{(k)} \frac{(e^t-1)^n}{n!} = \sum_{n=0}^{\infty} D_n^{(k)} \left\{ \sum_{m=n}^{\infty} S(m,n) \frac{t^m}{m!} \right\} = \sum_{m=0}^{\infty} \left\{ \sum_{n=0}^{m} D_n^{(k)} S(m,n) \right\} \frac{t^m}{m!} .$$

This gives

$$\sum_{m=0}^{\infty} B_m^{(k)} \frac{t^m}{m!} = \sum_{m=0}^{\infty} \left\{ \sum_{n=0}^{m} D_n^{(k)} S(m,n) \right\} \frac{t^m}{m!} .$$

Comparing coefficients we find

$$B_m^{(k)} = \sum_{n=0}^{m} D_n^{(k)} S(m,n) = \sum_{n=0}^{m} s(m+k,k) \binom{m+k}{k}^{-1} S(m,n)$$

as needed.

The following representation was recently published by Khattou et al. [20].

**Proposition 7.** *We have for* $0 \le n \le k-1$

(9) $$B_n^{(k)} = (-1)^n \, | s(k,k-n) | \binom{k-1}{n}^{-1} .$$

## References

**[1]** **Abdelmejid Bayad, Yoshinori Hamahata,** Polylogarithms and poly-Bernoulli polynomials, *Kyushu J. Math*., 65 (2011) 15-24.

**[2]** **Khristo N. Boyadzhiev**, Close encounters with the Stirling numbers of the second kind, *Math. Magazine,* 85 (4) (2012), 252-266.

**[3]** **Khristo N. Boyadzhiev**, Power series with binomial sums and asymptotic expansions, *Int. J. Math. Analysis,* Vol. 8 (2014), no. 28, 1389-1414.

**[4]** **Andei Z. Broder**, The r-Stirling numbers, *Discrete Math*., 49 (1984), 241-259.

**[5]** **Louis Comtet**, *Advanced Combinatorics*, Kluwer, 1974.

**[6]** **H. L. Garabedian,** A new formula for the Bernoulli numbers, Bull. Amer. Math. Soc., 46 (6) (1940), 531-533/

[7] **Henry W. Gould**, Explicit formulas for Bernoulli numbers, *Amer. Math. Monthly*, 79 (1) (1972), 44-51.

[8] **Ronald I. Graham, Donald E. Knuth, and Oren Patashnik**, *Concrete Mathematics*, (Addison-Wesley, New York) (1994)

[9] **Jesus Guillera, Jonathan Sondow**, Double integrals and infinite products for some classical constants via analytic continuations of Lerch's transcendent, *Ramanujan J*. 16 (2008) 247-270

[10] **B.-N. Guo, I. Mez˝o, and F. Qi**, An explicit formula for Bernoulli polynomials in terms of r-Stirling numbers of the second kind, *Rocky Mountain J. Math*., 46(6) (2016), 1919-1923.

[11] **Helmut Hasse,** Ein Summierungsverfahren fur die Riemannische $\zeta$ -Reihe, *Math. Z,* 32 (1930) 458-464.

[12] **Masanobu Kaneko**, Poly-Bernoulli numbers, *J. Theorie Nombres Bordeaux*, Vol. 9 (1997), 199-206.

[13] **Dae Sun Kim et al**., Higher-order Daehee numbers and polynomials, *Int. Journal Math. Anal*. 8(6) (2014), 273-283.

[14] **Takao Komatsu**, Poly-Cauchy numbers, *Kyushu J. Math*., Vol. 67 (2013), 143-153.

[15] **Takao Komatsu, Istvan Mezo**, Several explicit formulae of Cauchy polynomials in terms of r-Stirling numbers, *Acta Math. Hungar*.,148 (2) (2016), 522-529.

[16] **Donatella Merlini, Renzo Sprugnoli, M. Cecilia Verri**, The Cauchy numbers, *Discrete Math.,* 306 (2006), 1906-1920.

[17] **Antonio F. Neto**, A Note on a Theorem of Guo, Mezo, and Qi, *J. Integer Seq*., 19 (2016), 16.4.8.

[18] **Niels Nielsen**, Traité élémentaire des nombres de Bernoulli, Gauthier-Villars, Paris, 1923.

[19] **Pavel G. Todorov,** On the theory of the Bernoulli polynomials and numbers, *J. Math. Anal. Appl.* 104 (1984) 309-350.

[20] **Chouaib Khatttou, Abdelmejid Bayad, and Mohand Quamar Hernane,** New results on Bernoulli numbers of higher order, *Rocky Mountain J. Math*., 52 (2022), 153-170.